\documentclass[a4paper,DIV=12,11pt,parskip=half,notitlepage,abstracton]{scrartcl}
\usepackage[T1]{fontenc}
\usepackage[utf8]{inputenc}

\usepackage[pdfencoding=auto,
unicode=true,
plainpages=false,
pdfpagelabels=true]{hyperref}

\usepackage{graphicx}

\usepackage{booktabs}

\usepackage{siunitx}

\usepackage{wrapfig}

\usepackage{amssymb,amsmath}
\usepackage{stmaryrd}
\usepackage{mathtools}
\renewcommand{\d}{\mathrm d}
\renewcommand{\vec}{\boldsymbol}

\newcommand{\R}{\mathbb{R}}

\newcommand{\E}{\mathcal{E}}
\usepackage{booktabs}

\graphicspath{{graphics/}}

\usepackage{cleveref}

\usepackage{caption}
\usepackage{subcaption}

\begin{document}
\title{Numerical convergence of discrete extensions in a space-time finite element, fictitious domain method for the Navier--Stokes equations}

\author{Mathias Anselmann$^{\star,1}$, Markus Bause$^{\star}$\\
	{\small ${}^\star$ Helmut Schmidt University, Faculty of 
		Mechanical Engineering, Holstenhofweg 85,}\\ 
	{\small 22043 Hamburg, Germany}\\
}

\date{}
\maketitle

\footnotetext[1]{Corresponding author: anselmann@hsu-hh.de}

\vspace*{-8ex}
\begin{abstract}
A key ingredient of our fictitious domain, higher order space-time cut finite element (CutFEM) approach for solving the incompressible Navier--Stokes equations on evolving domains (cf.\ \cite{Bause2021}) is the extension of the physical solution from the time-dependent flow domain $\Omega_f^t$ to the entire, time-independent computational domain $\Omega$. The extension is defined implicitly and, simultaneously, aims at stabilizing the discrete solution in the case of unavoidable irregular small cuts. Here, the convergence properties of the scheme are studied numerically for variations of the combined extension and stabilization.  
\end{abstract}
\maketitle                   

\section{Mathematical problem and numerical scheme}
\label{Sec:Scheme}

\begin{figure}[!ht]
	\centering
	\captionsetup{width=.8\linewidth}
	\includegraphics[width=0.3\textwidth,keepaspectratio]{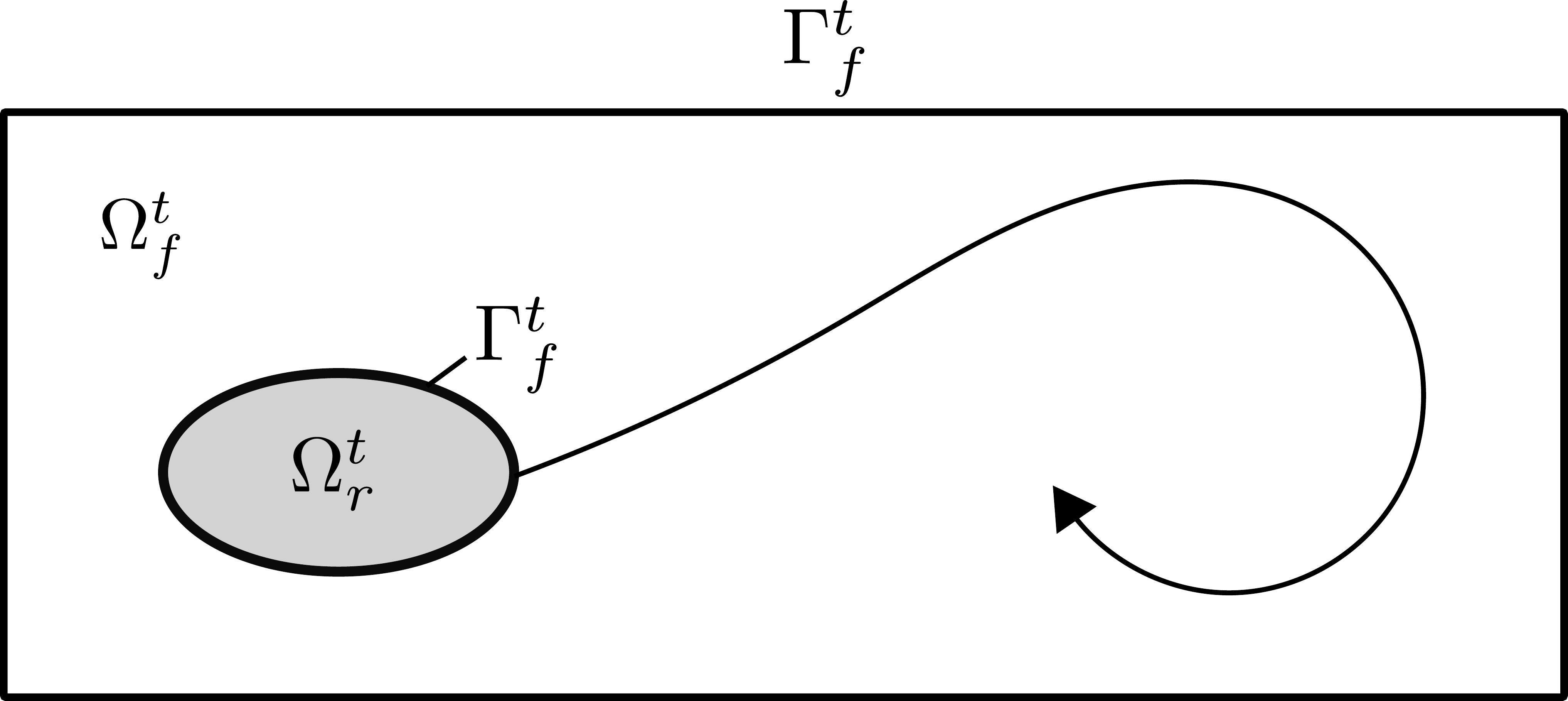}
	\caption{Prototype flow domain $\Omega_f^t$ with enclosed rigid body $\Omega_r^t$, computational background domain $\Omega$.}
	\label{fig:problem_setup}
\end{figure}

For an evolving flow domain $\Omega_f^t$, with $t\in [0,T]$, and $I:=(0,T]$ we consider the incompressible Navier--Stokes system (cf.\ Fig.\ \ref{fig:problem_setup})
\begin{subequations}
	\label{eq:navier_stokes}
	\begin{align*}
		\partial_t \vec{v} + (\vec{v} \cdot \nabla) \vec{v}  - \nu \Delta \vec{v} + \nabla p = \vec{f}\,, \quad 
		\nabla \cdot \vec{v} &= 0
		\hspace*{3ex} \text{in } \Omega_f^t  \times I\,,
	\end{align*}
\end{subequations}
that is equipped with the initial condition $\vec{v}(0) = \vec{v}_0$ in $\Omega_f^0$ and the Dirichlet boundary condition $\vec{v} = \vec{g}$ on $\Gamma_f^t  \times I$ for the time-dependent boundary $\Gamma_f^t$ of the flow domain $\Omega_f^t$.
For flow problems with inflow and outflow boundaries we refer to cf.\ \cite{Bause2021}.

Let $\mathcal T_h= \{K\}$ be a family of regular decompositions of the computational domain $\Omega=\Omega_f^t\cup \overline{\Omega_r^t}$. By $\vec V_h\times Q_h$ we denote an inf-sup stable pair of  time-independent bulk finite element spaces on $\Omega$ for the velocity and pressure.
For $\vec u_h, \vec \phi_h \in \vec V_h\times Q_h$, with $\vec u_h=(\vec v_h, p_h)$, $\vec \phi_h = (\vec \psi_h,\xi_h)$, we define $A_h: (\vec V_h\times Q_h)\times (\vec V_h\times Q_h)\rightarrow \R$ and $L_h: (\vec V_h\times Q_h)\rightarrow \R$ by
\begin{align*}
	\label{Def:Ah}
	A_h(\vec u_h,\vec \phi_h) & :=
	\begin{aligned}[t]
	&\langle (\vec v \cdot \vec \nabla) \vec v, \vec \psi \rangle_{\Omega_f^t}
	+ \nu \langle \nabla \vec v , \nabla \vec \psi  \rangle_{\Omega_f^t}
	-\langle p, \nabla \cdot \vec \psi \rangle_{\Omega_f^t}
	+ \langle \vec \nabla \cdot \vec v, \xi \rangle_{\Omega_f^t}
	\\
	&- \langle \nu \nabla \vec v_h \cdot \vec n - p_h \vec n, \vec \psi_h\rangle_{\Gamma_f^t}
	+ B_{\Gamma_f^t}(\vec v_h,\vec \phi_h)
	+ S_{F_h^t}(´\vec u_h,\vec \phi_h) \,,
	\end{aligned}
	\\
	L_h(\vec \phi_h;\vec f, \vec g)	&:=  \langle \vec f, \vec \psi_h\rangle_{\Omega_f^t} + B_{\Gamma_f^t}(\vec g,\vec \phi_h) \,.
\end{align*}
The bilinearform  $B_{\Gamma_f^t} : \vec H^{1/2}(\Gamma_f^t) \times (\vec V_h \times Q_h) \rightarrow \R$ enforces Dirichlet boundary conditions weakly by Nitsche's method,
\begin{equation*}
	\label{Def:B_GamD}
		B_{\Gamma_f^t}(\vec w,\vec \phi_h) : =   - \langle \vec w, \nu \nabla \vec \psi_h \cdot \vec n + \xi_h  \vec n  \rangle_{\Gamma_D^t}
		+ \gamma_1 \nu \langle h^{-1} \vec w , \vec 
		\psi_h  \rangle_{\Gamma_f^t}  + \gamma_2 \langle h^{-1} \vec w \cdot \vec n, \vec \psi_h \cdot \vec n 
		\rangle_{\Gamma_f^t}\,, 
\end{equation*}
for $\vec w \in \vec H^{1/2}(\Gamma_f^t)$ and $\vec \phi_h \in \vec V_h \times Q_h$. Here,
$\gamma_1>0$ and $\gamma_2> 0$ are numerical (tuning) parameters for the penalization. The linear form $S_{F_h^t}$ is introduced and analyzed for the Stokes problem in \cite{Wahl2020}. Its impact is twofold. On the one hand it extends the solution $\vec{u}_h$ from the fluid domain $\Omega_f^t$ to the rigid domain $\Omega_r^t$ such that the solution is defined on the whole computational domain $\Omega = \Omega_{f}^t \cup \overline{\Omega_{r}^t}$. On the other hand, $S_{F_h^t}$ stabilizes the solution in the case of small cut cell scenarios. For the definition of $S_{F_h^t}$, we define a stabilization zone $\Omega_{s}^t\subset \Omega $ such that the enclosure $\Omega_{r}^t \subset \Omega_{s}^t$ is satisfied, and a corresponding submesh $\mathcal T_{h,s}^{t}\subset \mathcal T_h$ of all cells that cover $\Omega_{s}^t$ completely. In Sec.\ 2 we will study numerically the effect of extending the stabilization into the fluid domain, i.e.\ of widening the domain $\Omega_{s}^t$.
The set of all faces that are common to two cells of $\mathcal T_{h,s}^{t}$ is denoted by $F_h^t$.
With numerical parameters $\gamma_v= \widetilde{\gamma}_v({\nu}^{-1} + \nu) h^{-2}>0$ and $\gamma_p= \widetilde{\gamma}_p{\nu}^{-1}>0$, the form $S_{F_h^t}$ is defined by 
\begin{align*}
	\small
	S_{F_h^t}(\vec u_h,\vec \phi_h)
	=
	\sum_{F \in F_h^t}
	\gamma_v \langle \E \vec{v}_{h|K_1} - \E \vec{v}_{h|K_2},
	\E \vec{\psi}_{h|K_1} - \E \vec{\psi}_{h|K_2}
	\rangle_{\omega_F}
	+ \gamma_p 
	\langle 
	\E p_{|K_1} - \E p_{|K_2},
	\E \xi_{|K_1} - \E \xi_{|K_2}
	\rangle_{\omega_F}\,.
\end{align*}
Here, $\E$ denotes the canonical patchwise extension of the discrete functions. We refere to \cite{Wahl2020,Bause2021} for its definition. 

Putting $A_h^s \coloneqq A_h + S_{F_h^t}$ and applying a discontinuous Galerkin time discretization with piecewise polynomials of order $k$ (cf.\ \cite{Bause2021}) leads to finding, in each subintervall $I_n=(t_{n-1},t_n)$, functions $(\vec v_{\tau,h},p_{\tau,h})\in \mathbb P_k(I_n;\vec V_{h}) \times \mathbb P_k(I_n;Q_h)$, such that 
\begin{align*}
	\int_{t_{n-1}}^{t_{n}} \langle\partial_t \vec v_{\tau,h}, \vec \psi_{\tau,h} \rangle_{\Omega_f^t} +
	A_h^s(\vec u_{\tau,h},\vec \phi_{\tau,h}) \d t + \langle [\vec v_{\tau,h}]_{t_{n-1}}, \vec \psi _{\tau,h}(t_{n-1}^+) \rangle_{\Omega}  = \int_{t_{n-1}}^{t_{n}} 
	L_h(\vec \phi_{\tau,h}; \vec f, \vec g) \d t 
\end{align*}
for all $\vec \phi_{\tau,h}\in \mathbb P_k (I;\vec V_h)\times \mathbb P_k(I_n;Q_h)$, with $\vec \phi_{\tau,h}= (\vec \psi_{\tau,h},\xi_{\tau,h})$.

\section{Numerical convergence study}

\vspace*{-1ex}
\begin{figure}[!htb]
	\centering
	\subcaptionbox{\centering Problem setting. \label{fig:experiment_setup}}
	[0.32\columnwidth]
	{\includegraphics[width=0.185\textwidth,keepaspectratio]
		{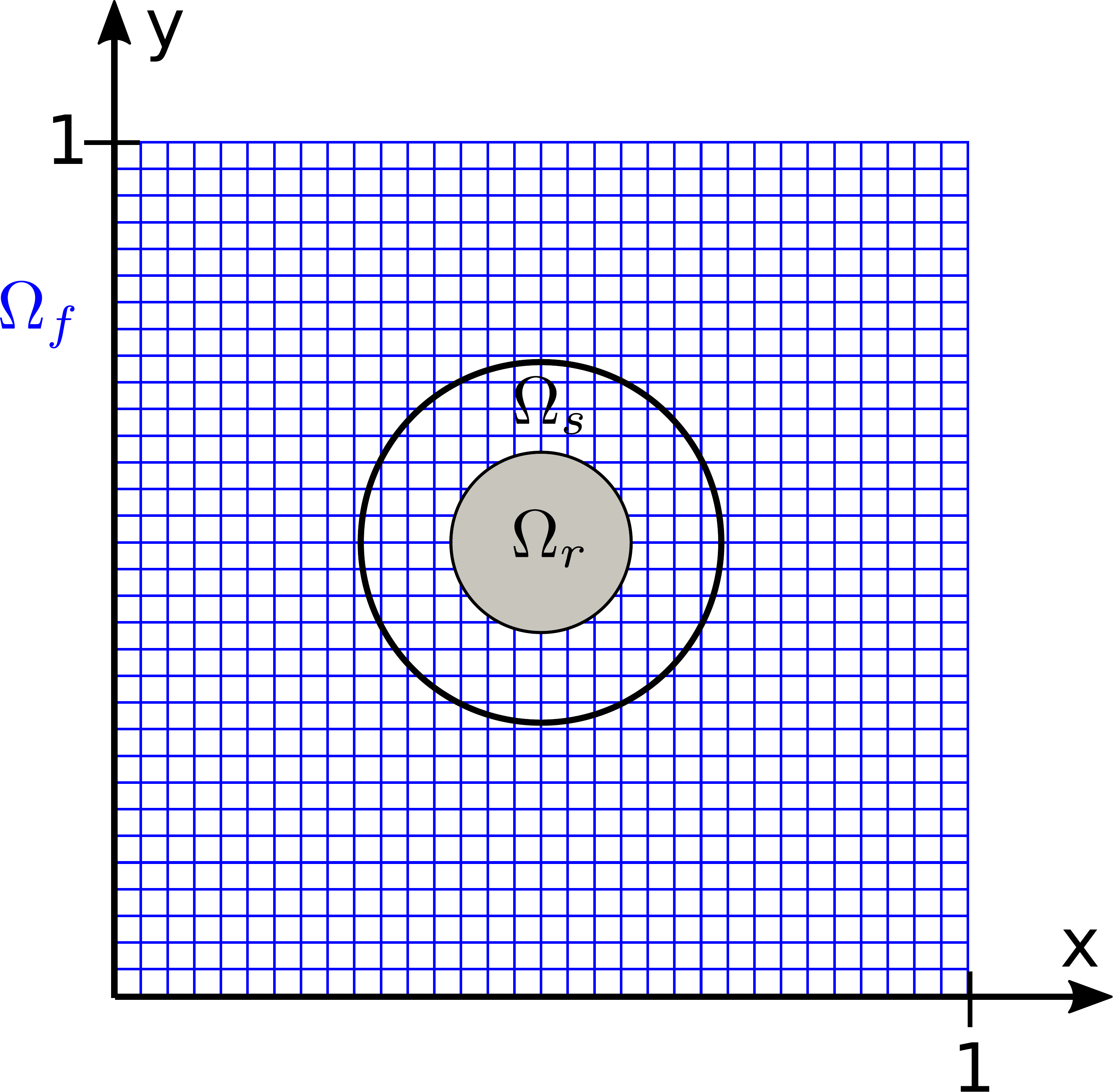} }
	\subcaptionbox{\centering Submesh $\mathcal T_{h,s}^{t}$ of stabilization \label{fig:stabilization_cells}}
	[0.32\columnwidth]
	{\includegraphics[width=0.175\textwidth,keepaspectratio]
		{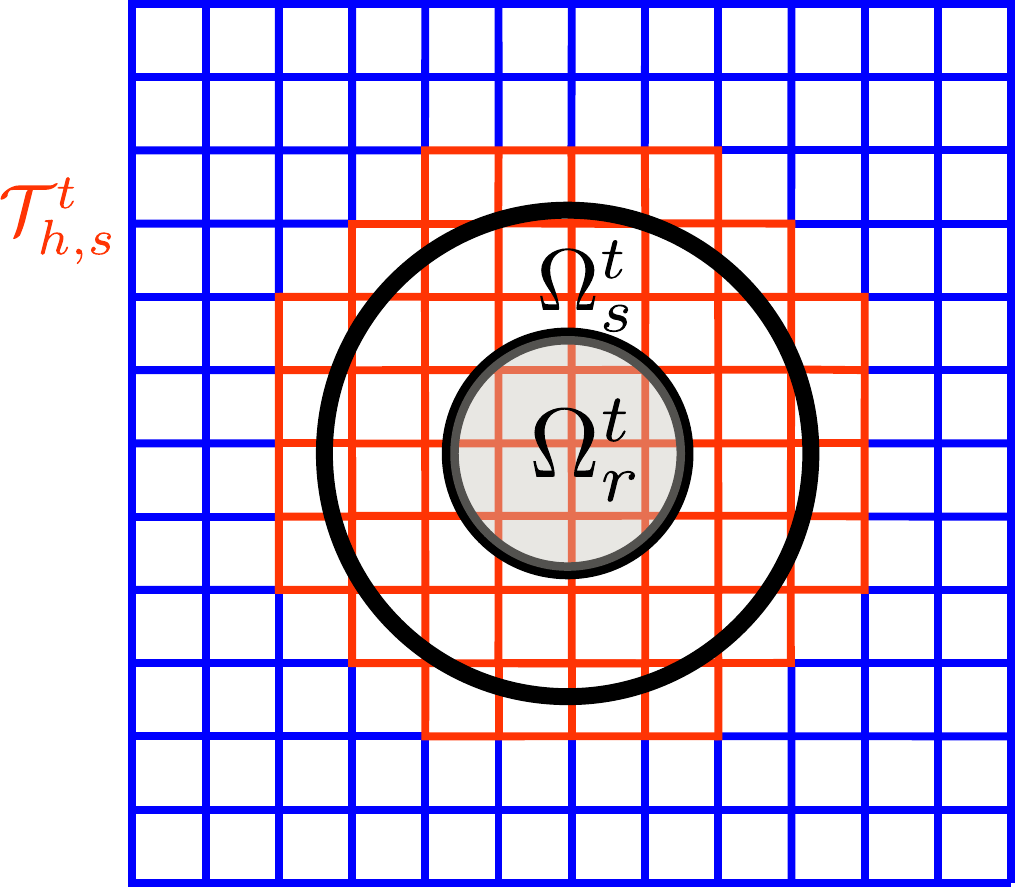}} 
	\subcaptionbox{\centering Pressure field at $t = 1$ for $h = {h_0}/{2^4}$. \label{fig:pressure_field}}
	[0.32\columnwidth]
	{\includegraphics[width=0.175\textwidth,keepaspectratio]{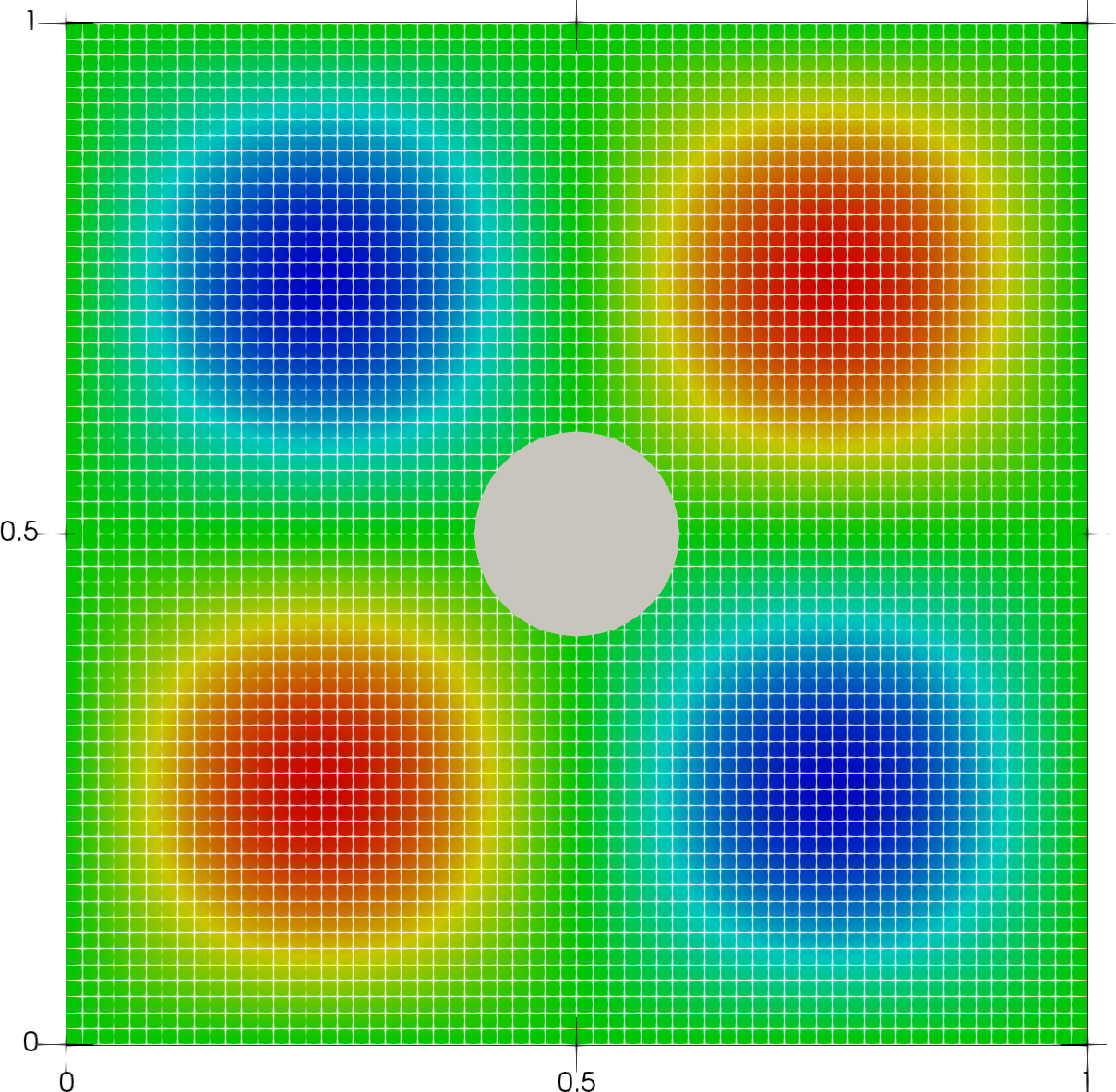}}
	\caption{Problem setting, domain $\Omega_{s}^t$ and submesh $\mathcal T_{h,s}^{t}$ of stabilization and computed solution profile.}\label{fig:numerical_experiments}
\end{figure}

For the sake of implementational simplicity, we consider the problem setting of Fig.~\ref{fig:experiment_setup} with a time-independent domain $\Omega_r$. Our computational study investigates the impact of the width of the domain of stabilization $\Omega_{s}^t$ on the convergence of the scheme proposed in Sec.~\ref{Sec:Scheme}. Precisely, we aim to investigate if a wider overlapping of the fluid domain $\Omega_f^t$ by $\Omega_s^t$ is required to ensure convergence of optimal order. We put $\Omega \times I = (0,1)^2 \times (0,1]$. The midpoint of the circular rigid body with radius $r_{\Omega_r} = 0.1$ is located in the center of $\Omega$. We prescribe the initial value $\vec{v}_0$ and right-hand side function $\vec{f}$ on $\Omega_f \times I$ in such a way, that the solution of the Navier--Stokes system on $\Omega_{f}$ is given by
$\vec{v}_e(\vec{x},t)
\coloneqq
[ \cos(x_2 \pi) \cdot \sin(t) \cdot \sin(x_1 \pi)^2  \cdot \sin(x_2 \pi),$
$- \cos(x_1 \pi) \cdot \sin(t) \cdot \sin(x_2 \pi)^2  \cdot \sin(x_1 \pi) ]^\top$
and
$p_e(\vec{x},t)
\coloneqq
\cos(x_2 \pi) \cdot \sin(t) \cdot \sin(x_1 \pi) \cdot \cos(x_1 \pi) \cdot \sin(x_2 \pi).$
On the inner fluid boundary we prescribe the condition $\vec v= \vec{v}_e$, on the outer boundary we prescribe a homogeneous Dirichlet condition. Table~\ref{tab:numerical_results_1} shows the computed errors and experimental orders of convergence (EOC) for two different combinations of space-time finite elements, based on the Taylor--Hood family $(H_h^r)^2\times H_h^{r-1}$, $r\geq 2$ in space, and two different choices of the radius $r_{\Omega_s}$ of the domain $\Omega_s$ where the stabilization $S_{F_h^t}$ is applied, precisely $r_{\Omega_s} = 1\cdot r_{\Omega_r}$ (\textit{bottom}) and $r_{\Omega_s} = 2\cdot r_{\Omega_r}$ (\textit{top}).
\sisetup{scientific-notation = true,
		 round-mode=places,
		 round-precision=3,
		 output-exponent-marker=\ensuremath{\mathrm{e}},
		 table-figures-integer=1, 
		 table-figures-decimal=3, 
		 table-figures-exponent=1, 
		 table-sign-mantissa = false, 
		 table-sign-exponent = true, 
	 	 table-number-alignment=center} 
\begin{table}[!h]
	\caption{Errors and experimental order of convergence for varying radius $r_{\Omega_{s}}$ of the domain $\Omega_{s}$ for $\tau_0 = 1.0$ and $h_0 = {1}/({2 \sqrt{2}})$.}
	\centering
	\small
	\begin{tabular}{c@{\,\,\,\,}c  S@{\,}c  S@{\,}c  S@{\,}c  S@{\,}c}
		\toprule
		{$\tau$} & {$h$} &
		{ $\| e^{\;\vec{v}}  \|_{L^2(L^2)} $ } & {EOC} &
		{ $\| e^{\;p}  \|_{L^2(L^2)} $ } & {EOC} &  
		{ $\| e^{\;\vec{v}}  \|_{L^2(L^2)} $ } & {EOC} &
		{ $\| e^{\;p}  \|_{L^2(L^2)} $ } & {EOC} \\
		\cmidrule(r){1-2}
		\cmidrule(lr){3-4}
		\cmidrule(lr){5-6}
		\cmidrule(l){7-8}
		\cmidrule(l){9-10}
		$\tau_0/2^0$ & $h_0/2^0$ & 1.514e-02 & {--} & 3.719e-02 & {--} & 1.581e-03 & {--} & 8.851e-03 & {--}  \\ 
		$\tau_0/2^1$ & $h_0/2^1$ & 3.764e-03 & 2.01 & 9.213e-03 & 2.01 & 2.234e-04 & 2.82 & 1.174e-03 & 2.92  \\
		$\tau_0/2^2$ & $h_0/2^2$ & 9.100e-04 & 2.05 & 2.156e-03 & 2.09 & 2.972e-05 & 2.91 & 1.531e-04 & 2.94  \\
		$\tau_0/2^3$ & $h_0/2^3$ & 2.340e-04 & 1.96 & 5.524e-04 & 1.96 & 3.673e-06 & 3.02 & 1.962e-05 & 2.96  \\
		$\tau_0/2^4$ & $h_0/2^4$ & 5.792e-05 & 2.01 & 1.381e-04 & 2.00 & 4.666e-07 & 2.98 & 2.479e-06 & 2.98  \\
		$\tau_0/2^5$ & $h_0/2^5$ & 1.448e-05 & 2.00 & 3.477e-05 & 1.99 & 5.833e-08 & 3.00 & 3.120e-07 & 2.99  \\
		\cmidrule(r){1-2}
		\cmidrule(lr){3-6}
		\cmidrule(lr){7-10}
		\multicolumn{2}{c}{$\text{elements}_{|I_n}$} &
		\multicolumn{4}{c}{$(\mathbb P_1(I_n;H_h^2))^2 \times \mathbb P_1(I_n;H_h^1)$, $r_{\Omega_s} = 2 \cdot r_{\Omega_r}$} &
		\multicolumn{4}{c}{$(\mathbb P_2(I_n;H_h^3))^2 \times \mathbb P_2(I_n;H_h^2)$, $r_{\Omega_s} = 2 \cdot r_{\Omega_r}$} \\
		\bottomrule
		\addlinespace[2ex]
		\toprule
		{$\tau$} & {$h$} &
		{ $\| e^{\;\vec{v}}  \|_{L^2(L^2)} $ } & {EOC} &
		{ $\| e^{\;p}  \|_{L^2(L^2)} $ } & {EOC} &  
		{ $\| e^{\;\vec{v}}  \|_{L^2(L^2)} $ } & {EOC} &
		{ $\| e^{\;p}  \|_{L^2(L^2)} $ } & {EOC} \\
		\cmidrule(r){1-2}
		\cmidrule(lr){3-4}
		\cmidrule(lr){5-6}
		\cmidrule(l){7-8}
		\cmidrule(l){9-10}
		$\tau_0/2^0$ & $h_0/2^0$ & 1.540e-02 & {--} & 3.045e-02 & {--} & 1.552e-03 & {--} & 8.794e-03 & {--}  \\ 
		$\tau_0/2^1$ & $h_0/2^1$ & 5.762e-03 & 1.42 & 9.503e-03 & 1.68 & 2.191e-04 & 2.83 & 1.925e-03 & 2.19  \\
		$\tau_0/2^2$ & $h_0/2^2$ & 1.812e-03 & 1.67 & 2.248e-03 & 2.08 & 2.194e-05 & 3.32 & 1.743e-04 & 3.46  \\
		$\tau_0/2^3$ & $h_0/2^3$ & 5.004e-04 & 1.86 & 5.580e-04 & 2.01 & 5.293e-06 & 2.05 & 4.302e-05 & 2.03  \\
		$\tau_0/2^4$ & $h_0/2^4$ & 1.316e-04 & 1.93 & 1.444e-04 & 1.95 & 7.240e-07 & 2.87 & 5.723e-06 & 2.91  \\
		$\tau_0/2^5$ & $h_0/2^5$ & 3.290e-05 & 2.00 & 3.636e-05 & 1.99 & 9.240e-08 & 2.97 & 7.203e-07 & 2.99  \\
		\cmidrule(r){1-2}
		\cmidrule(lr){3-6}
		\cmidrule(lr){7-10}
		\multicolumn{2}{c}{$\text{elements}_{|I_n}$} &
		\multicolumn{4}{c}{$(\mathbb P_1(I_n;H_h^2))^2 \times \mathbb P_1(I_n;H_h^1)$, $r_{\Omega_s} = 1 \cdot r_{\Omega_r}$} &
		\multicolumn{4}{c}{$(\mathbb P_2(I_n;H_h^3))^2 \times \mathbb P_2(I_n;H_h^2)$, $r_{\Omega_s} = 1 \cdot r_{\Omega_r}$} \\
		\bottomrule
	\end{tabular}
	\label{tab:numerical_results_1}
\end{table}

Even though the absolute errors are smaller for $r_{\Omega_s} = 2\cdot r_{\Omega_r}$ (\textit{top}), our results indicate that the EOC does not depend on the diameter of the region where the stabilization $S_{F_h^t}$ is applied. We note that in our experiments cut cells with non-convex boundary segments due to intersections with $\Omega_r^t$ occur, cf.\ Fig.~\ref{fig:numerical_experiments}. The convergence does not suffer from these cells.


\begin{thebibliography}{1}

	\bibitem{Bause2021} M.\ Anselmann, M.\ Bause:
	CutFEM and ghost stabilization techniques for higher order space-time discretizations of the Navier--Stokes equations, submitted (2021),
	\href{https://arxiv.org/abs/2103.16249}{arxiv:2103.16249}.

    \bibitem{Bause2020} M.\ Anselmann, M.\ Bause:
    Higher order Galerkin–collocation time discretization with Nitsche’s method for the Navier--Stokes equations, Math.\ Comp.\ Simul. (2020),  \href{https://doi.org/10.1016/j.matcom.2020.10.027}{DOI:\,10.1016/j.matcom.2020.10.027}.
    
    \bibitem{Wahl2020}H.\ von\ Wahl, T.\ Richter, and C.\ Lehrenfeld:
    An unfitted Eulerian finite elementmethod for the time-dependent Stokes problem on moving domains,
    preprint (2020), \href{https://arxiv.org/abs/2002.02352}{arXiv:2002.02352}.
        
\end{thebibliography}
\end{document}